\documentclass{article}

\usepackage{amsmath,amssymb,amsfonts,amsthm,latexsym}
\usepackage{hyperref, enumerate}
\usepackage{tikz,tikz-network}
\usepackage{rotating}
\usepackage{ifthen}
\usepackage{geometry}
\usepackage{float}

\usetikzlibrary{calc, math}
\usetikzlibrary{decorations.markings}

\newtheorem{lemma}{\bf Lemma}
\newtheorem{remark}{\bf Remark}

\newtheorem{theorem}{\bf Theorem}

\newtheorem*{definition*}{\bf Definition}

\newtheorem*{problem*}{\bf Problem}

\newtheorem*{comment*}{\bf Comment}

\allowdisplaybreaks

\title{Resistance, oddness and colouring defect of snarks}

\author{Imran Allie}

\begin{document}
	
\maketitle
	
\begin{abstract}
	Let $G$ be a bridgeless cubic graph. The \textit{resistance} of $G$, denoted $r(G)$, is the minimum number of edges which can be removed from $G$ in order to render 3-edge-colourability. The \textit{oddness} of $G$, denoted $\omega(G)$, is the minimum number of odd components in a 2-factor of $G$. The \textit{colouring defect} of $G$ (or simply, the \textit{defect} of $G$), denoted $\mu_3(G)$, is the minimum number of edges not contained in any set of three perfect matchings of $G$. These three parameters are regarded as measurements of uncolourability of snarks, partly because any one of these parameters equal zero if and only if $G$ is 3-edge-colourable. It is also known that $r(G) \geq \omega(G)$ and that $\mu_3(G) \geq \frac{3}{2}\omega(G)$ \cite{fiol,jinsteffen}. We have shown that the ratio of oddness to resistance can be arbitrarily large for non-trivial snarks \cite{allie1}. It has also been shown that the ratio of the defect to oddness can be arbitrarily large for non-trivial snarks, although this result was only shown for graphs with oddness equal to 2 \cite{karabasetal}. In the same paper, the question was posed whether there exists non-trivial snarks for given resistance $r$ or given oddness $\omega$, and arbitrarily large defect. In this paper, we prove a stronger result: For any positive integers $r \geq 2$, even $\omega \geq r$, and $d \geq \frac{3}{2}\omega$, there exists a non-trivial snark $G$ with $r(G)=r$, $\omega(G)=\omega$ and $\mu_3(G) \geq d$.  
\end{abstract}

	
\section{Introduction}{\label{Introduction}}

\textit{Snarks} are bridgeless cubic graphs whose edges cannot be properly coloured with three colours. The significance of snarks derives mainly from the fact that it is the only class of graphs which may possibly contain counterexamples to many major open conjectures in graph theory. The cycle double cover conjecture, Tutte's $5$-flow conjecture, and Fulkerson's conjecture, are some examples of these major open conjectures. In particular, these potential counter-examples have also been shown, in some cases, to necessarily have girth greater than 4 and cyclic connectivity greater than 3. In this paper we thus continue the convention of referring to snarks with girth less than 5 or cyclic connectivity less than 4 as \textit{trivial}. These snarks are also considered trivial as they are easily reducible to smaller snarks by well-defined reductions. We do however believe that the idea of what constitutes triviality in snarks should be considered more deeply. See \cite{allie2} for further ideas on this discussion. In contemporary times, researchers have endeavoured to understand snarks by defining and investigating parameters which `measure', in various ways, how far a snark is from being 3-edge-colourable. Such parameters have been termed \textit{measures of uncolourability}. In this paper, we continue this endeavour by focusing on the parameters resistance, oddness and colouring defect. For related results on these parameters, see \cite{allie1,allie3,allie4,jinsteffen,karabasetal,kochol,steffen1}. See as well \cite{fiol} for a recent comprehensive on measurements of uncolourability.  

The \textit{(vertex) resistance} of $G$, denoted ($r_v(G)$) $r(G)$, is the minimum number of (vertices) edges which can be removed from $G$ in order to render 3-edge-colourability. (It is well-known that $r(G)=r_v(G)$ for any cubic $G$, and these two parameters could be used interchangeably). The \textit{oddness} of $G$, denoted $\omega(G)$, is the minimum number of odd components in a 2-factor of $G$. These are examples of measures which have been extensively studied, specifically in relation to each other. In fact, in \cite{fiol} Fiol et al conjectured that $\omega(G) \leq 2r(G)$ for any bridgeless cubic graph. In \cite{allie1}, we disproved this conjecture by presenting a class of non-trivial snarks in which the ratio of oddness to resistance can be arbitrarily large. Specifically, we proved that there exists no $k$ such that $\omega(G) \leq kr(G)$ for any snark $G$. 

Another measure is the colouring defect of $G$. The \textit{colouring defect} of $G$ (or simply, the \textit{defect} of $G$), denoted $\mu_3(G)$, is the minimum number of edges not contained in any set of three perfect matchings of $G$. This measure was introduced by Steffen in \cite{steffen1}. Steffen proved, among other things, that any snark has defect at least 3. Another notable result states that the defect of a snark is at least as large as one half of its girth. Since there exists snarks of arbitrarily large girth \cite{karabasetal,kochol}, there exist snarks of arbitrarily large defect. In \cite{jinsteffen}, Jin and Steffen proved that $\mu_3(G) \geq \frac{3}{2}\omega(G)$. Jin and Steffen additionally considered the ratio difference of the two parameters, and to this end proved that for any given oddness $\omega$ and any $d \geq \frac{3}{2}\omega$ there exists a bridgeless cubic graph with oddness $\omega$ and defect at least $d$. However, their construction produces snarks with cyclic connectivity not exceeding 3 (thus trivial snarks). In \cite{karabasetal}, the authors nearly completely generalise this result by Jin and Steffen to cyclically 5-edge-connected snarks (thus non-trivial snarks). The generalisation is only lacking in completeness in that the snarks only have oddness equal to 2 in each instance. The authors mention that their graphs are not simply adjustable to have arbitrary oddness, together with the same desired property of non-triviality and arbitrarily large girth. The authors then pose the question about whether such graphs exist, and suggest an answer in the affirmative. That is, whether there exists non-trivial snarks with any given possible oddness and arbitrarily large defect, or the slightly weaker version which asks whether there exists non-trivial snarks with any given possible resistance and arbitrarily large defect.

In this paper, we completely generalise the result by Jin and Steffen for non-trivial snarks by answering the question posed in \cite{karabasetal}. In fact, we prove a stronger result which involves resistance, oddness and defect. We prove that for any integers $r \geq 2$, even $\omega \geq r$, and $d \geq \frac{3}{2}\omega$, there exists a cyclically 4-edge-connected snark $G$ with girth 5 (thus non-trivial) such that $r(G)=r$, $\omega(G)=\omega$ and $\mu_3(G) \geq d$. We do this by presenting an infinite class of graphs in which the resistance is easily adjustable, oddness is easily adjustable whilst keeping resistance constant, and defect is easily ensured to be arbitrarily large by increasing the girth of just a conflicting subgraph of the snark. This construction relies heavily on the methods of snark construction presented in \cite{allie1} and \cite{karabasetal}.

\section{Preliminaries}{\label{preliminaries}}

A \textit{semi-graph} $G$ is a pair $G=(V,E)$ which consists of a set of vertices $V = V(G)$ and a set $E = E(G) \subseteq \mathcal{P}_2(V) \cup \mathcal{P}_1(V)$ consisting of edges and semi-edges; here $\mathcal{P}_i(A)$ denotes the set of all $i$-element subsets of a set $A$. In $E(G)$, the 2-element sets are called \textit{edges} (as expected) while the 1-element sets are called \textit{semi-edges}. Note that if $E$ contains no elements from $\mathcal{P}_1(V)$, then $G$ is simply a \textit{graph}.

We denote the edge $\{u, v\}$ as $uv$ and the semi-edge $\{u\}$ as $(u)$. Furthermore, we define the \textit{join} between two semi-edges $(u)$ and $(v)$ as the removal of semi-edges $(u)$ and $(v)$, and the addition of the edge $uv$. A semi-edge $(u)$ and a vertex $v$ may also join to form an edge $uv$, with semi-edge $(u)$ being removed. The degree of a vertex $v$ in a semi-graph $G$ is defined as the combined total number of edges and semi-edges incident with $v$. Thus a cubic semi-graph is a semi-graph with each vertex having degree 3. Cycles and paths are defined analogously in semi-graphs as in graphs and may not contain semi-edges except that a \textit{$k$-semi-path} contains $k$ semi-edges (obviously, it is only possible that $k=1$ or $k=2$ and any semi-edge in a $k$-semi-path would necessarily be terminal). The length of a semi-path is the combined total number of edges and semi-edges. The \textit{semi-girth} of $G$ is the minimum order of a cycle or 2-semi-path in $G$. Essentially, semi-edges behave like edges except that they are associated with one vertex instead of two. We say that a semi-graph $G$ contains a semi-subgraph $G'$ if $V(G') \subseteq V(G)$, $uv \in E(G')$ implies that $uv \in E(G)$, semi-edge $(u) \in E(G')$ implies semi-edge $(u) \in E(G')$ or there is an edge $uv \in E(G)$, and for every vertex $u \in V(G')$ the degree of $u$ in $G$ is greater than or equal to the degree of $u$ in $G'$ (as a quick example, the semi-graph with one vertex and one semi-edge is a semi-subgraph of $K_2$). We note that all semi-graphs considered in this paper are finite, loopless, have no parallel edges, and each vertex is associated with at most one semi-edge.

Let $G=(V,E)$ be a graph. A \textit{$k$-edge-colouring}, $f$, of $G$ is a mapping from the set of edges of $G$ to a set of $k$ colours. That is,  $f : E \longrightarrow \{1,\dots,k\}$. A $k$-edge-colouring of $G$ is \textit{proper} if no two adjacent elements in $E$ are mapped to the same colour. By Vizing's theorem \cite{vizing1,vizing2}, the smallest number $k$ such that a simple graph $G$ admits a proper $k$-edge-colouring is either $\Delta(G)$ or $\Delta(G) + 1$, where $\Delta$ is the maximum degree of any vertex in $G$. It is easy to see that Vizing's theorem applies for semi-graphs as well. A vertex $v$ is \textit{conflicting} with regard to an edge-colouring $f$ if more than one of the edges incident to $v$ are mapped to the same colour. The vertex resistance of cubic $G$ can thus also be thought of as the minimum number of conflicting vertices of a 3-edge-colouring of $G$.

The definition of a $k$-factor of a graph can also naturally be extended to a semi-graph. That is, a \textit{$k$-factor of a semi-graph $G$} is a $k$-regular semi-subgraph of $G$. Petersen's theorem is a well-known result which states that any bridgeless cubic graph $G$ contains a 1-factor \cite{petersen}. Moreover, Schonberger extended this result by proving that for any $e \in E(G)$ there exists a 1-factor of $G$ which contains $e$ \cite{schonberger}. We confirm the generalisation of this result to cubic semi-graphs in Lemma \ref{lem_cubicsemigraphhas1factor}. Lemma \ref{lem_cubicsemigraphhas1factor} is used implicitly throughout the paper.

\begin{lemma} \label{lem_cubicsemigraphhas1factor}
	Let $G$ be a bridgeless cubic semi-graph. Then for any $e \in E(G)$, $G$ has a 1-factor which contains $e$.
	\begin{proof}
		Let $e \in G$. By joining the semi-edges of $G$ (adding a vertex if necessary), we render a bridgeless cubic graph $G'$. $G'$ then has a 1-factor $M'$ which contains $e$ (if $e$ was a semi-edge in $G$, then it is now joined with another semi-edge or vertex and is now an edge in $G'$). We now reverse the process of joining edges (and adding the vertex) and define a set $M \subset E(G)$ as follows: if $e \in M'$ and $e$ was an edge in $G$, then $e \in M$; if $e \in M'$ and $e$ was formed by joining a semi-edge $e_1$ to an added vertex in $G$, then $e_1$ is in $M$; if $e \in M'$ and $e$ was formed by the joining of two semi-edges $e_1$ and $e_2$ in $G$, then $e_1$ and $e_2$ is in $M$. $M$ is then a 1-factor of $G$ which contains $e$.
	\end{proof}
\end{lemma}

The rest of this section presents and generalises the theory of Jin and Steffen in \cite{jinsteffen} from cubic graphs to cubic semi-graphs, and is useful in proving the main result of this paper. Let $G$ be a bridgeless cubic graph and let $\mathcal{M}=\{M_1,M_2,M_3\}$ be a set of three 1-factors of $G$. Let $E_i[\mathcal{M}]$ be the set of edges of $G$ which are in exactly $i$ of the three 1-factors in $\mathcal{M}$. Let $M = E_2 \cup E_3$, $U = E_0$ and $|U| = k$. (Clearly, the minimum possible $k$ coincides with the definition of $\mu_3(G)$. If $k = \mu_3(G)$ then $\mathcal{M}$ is an \textit{opimal} array of 1-factors). The \textit{$k$-core}, or simply, the \textit{core of $G$ with respect to $\mathcal{M}$}, is the subgraph $G_c$ of $G$ which is induced by $M \cup U$. A core is \textit{proper} if $G_c$ is a strict subset of $G$. Steffen proved various properties of cores in \cite{steffen1}. It was shown that every bridgeless cubic graph $G$ has a proper core. Obviously  then, every core of $G$ such that $k = \mu_3(G)$ is proper. 

\begin{lemma} \label{lem_cubicsemigraphhascore}
	Let $G$ be a cubic semi-graph. Then $G$ has a proper core with respect to some set of three 1-factors. 
	\begin{proof}
		Let $G'$ be a cubic graph obtained from joining semi-edges of $G$, adding a vertex if necessary. Let $\mathcal{M'}=\{M'_1,M'_2,M'_3\}$ be a set of three 1-factors of $G'$ such that the core of $G'$ with respect to $\mathcal{M'}$ is proper. Let $\mathcal{M}=\{M_1,M_2,M_3\}$ be a set of three 1-factors of $G$ inherited from the set $\mathcal{M'}$ of 1-factors of $G'$. Since the core of $G'$ with respect to $\mathcal{M'}$ is proper, it is clear that the core of $G$ with respect to $\mathcal{M}$ is also proper.   
	\end{proof}
\end{lemma}

Since a cubic semi-graph admits a 1-factor, we note that we can define the $defect$ of a cubic semi-graph just as it is defined for a cubic graph. That is, the minimum number of edges not contained in any set of three perfect matchings. Also, for a cubic semi-graph $G$, it is easy to see that $\mu_3(G) = 0$ if and only if $G$ is 3-edge-colourable. 

Steffen also showed that every component of a core of a cubic graph $G$ is either a cycle or a subdivision of a cubic graph. We generalise this property to semi-graphs as well.

\begin{lemma} \label{lem_prelim1}
	Let $G$ be a cubic semi-graph. Each component of the core of $G$ is a subdivision of a cubic semi-graph, an even cycle, or a 2-semi-path.    
	\begin{proof}		
		Let $\mathcal{M}=\{M_1,M_2,M_3\}$ be a set of three 1-factors of $G$. Let $M$ and $U$ be as previously defined with respect to $\mathcal{M}$. Let $G_c$ be the core of $G$ with respect to $\mathcal{M}$ and let $K_c$ be a component of $G_c$. If $K_c$ has no trivalent vertices, then $E(K_c) \cap M$ is a 1-factor of $K_c$ and hence, $K_c$ is an even cycle or a 2-semi-path whose edges and semi-edges are in $M$ and $U$ alternately. If $K_c$ contains trivalent vertices, then it is a subdivision of a cubic semi-graph.
	\end{proof}
\end{lemma}

Since the core of a cubic graph is a cycle or subdivision of a cubic graph, at least one edge incident to every 2-valent vertex in any $\mu_3(G)$-core of a cubic graph $G$ is uncovered by any of the 1-factors, and at least two edges incident to every 3-valent vertex in any $\mu_3(G)$-core of $G$ is uncovered by any of the 1-factors, it follows that $\mu_3(G) \geq \lceil \frac{\text{girth}(G)}{2} \rceil$. By the same logic and Lemma \ref{lem_prelim1}, it is then easy to see that we have the following analogous result for a cubic semi-graph $G$: $\mu_3(G) \geq \lceil \frac{\text{semi-girth}(G)}{2} \rceil$. In fact, we prove a stronger version relating to conflicting subgraphs. A \textit{conflicting subgraph} of a cubic semi-graph $G$ is a semi-subgraph which is itself not 3-edge-colourable. Conflicting subgraphs were introduced in \cite{fiol} (called conflicting zones) and minimal conflicting subgraphs of snarks were introduced and studied in \cite{allie2}. The following Lemma is at the crux of the proof of the main result.

\begin{lemma} \label{lem_semigirth}
	Let $G$ be a snark which contains a conflicting cubic semi-subgraph $H$. Then $\mu_3(G) \geq \lceil \frac{\text{semigirth}(H)}{2} \rceil$. 
	\begin{proof}
	We know that $\mu_3(H) \geq \lceil \frac{\text{semi-girth}(H)}{2} \rceil$. Assume that $\lceil \frac{\text{semigirth}(H)}{2} \rceil > \mu_3(G)$. Then $\mu_3(H) > \mu_3(G)$. Let $\mathcal{M}$ be an optimal array of 1-factors of $G$. Let $\mathcal{M'}$ be the array of $H$ inherited from $\mathcal{M}$. Let $M'$ and $M$ be the set of edges in $\mathcal{M'}$ and $\mathcal{M}$ not contained in any of the 1-factors, respectively. Clearly, $\mu_3(G) = |M| \geq |M'| \geq \mu_3(H)$ contradicting that $\mu_3(H) > \mu_3(G)$. Therefore, $\mu_3(G) \geq \lceil \frac{\text{semigirth}(H)}{2} \rceil$.
	\end{proof}
\end{lemma}

\section{Main result}{\label{mainresult}}

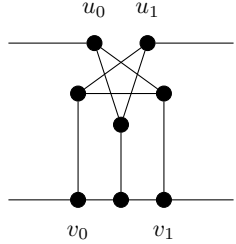
\begin{figure}[H] 
	\begin{center}
		\begin{tikzpicture}[every node/.style={draw,shape=circle,fill=black,text=white,scale=0.6},scale=0.4]
			\foreach \i [count=\ii from 0] in {54,126,198,270,342}{
				\path (\i:15mm) node (p\ii) {};
			}
			\foreach \x in {0,...,4}{
				\tikzmath{
					integer \y;
					\y = mod(\x+2,5);
				}
				\draw (p\x) -- (p\y);
			}
			\node [fill=white,draw=none] (b0) at (-4,-4) {};
			\node [fill=white,draw=none] (b4) at (4,-4) {};
			\draw (b0) -- (b4);
			\node (b1) at ($(b0)!(p2)!(b4)$) {};
			\node (b2) at ($(b0)!(p3)!(b4)$) {};
			\node (b3) at ($(b0)!(p4)!(b4)$) {};
			\draw (p2) -- (b1);
			\draw (p3) -- (b2);
			\draw (p4) -- (b3);
			\path let \p1 = (p1) in node [fill=white,draw=none] (a0) at (-4,\y1) {};
			\path let \p1 = (p1) in node [fill=white,draw=none] (a1) at (4,\y1) {};
			\draw (p1) -- (a0);
			\draw (p0) -- (a1);
			\path [late options={name=p0,label={[color=black,scale=1.5]:$u_1$}}];
			\path [late options={name=p1,label={[color=black,scale=1.5]90:$u_0$}}];
			\path [late options={name=b1,label={[color=black,scale=1.5]270:$v_0$}}];
			\path [late options={name=b3,label={[color=black,scale=1.5]270:$v_1$}}];
		\end{tikzpicture}
	\end{center}
	\caption{Cubic semi-graph $X$. The Petersen graph with two adjacent vertices removed.}\label{x}
\end{figure}

As mentioned, the proof of our main result relies heavily on constructions of graphs presented in \cite{allie1} and \cite{karabasetal}. We present these constructions and pertinent information regarding them before we formally begin proving the main result. Figure \ref{x}, Lemma \ref{lemmaX} and Remark \ref{remarkX} were presented in \cite{allie1}.

\begin{lemma} \label{lemmaX} 
	Consider the semi-graph $X$ in Figure \ref{x}. The following statements are true. 
	\begin{enumerate} 
		\item[(i)] There exists Hamiltonian paths in $X$ with and only with end-vertices $u_0$ and $v_0$, or $v_1$ and $u_1$. 
		\item[(ii)]  $X$ is Hamiltonian. 
		\item[(iii)] The girth of $X$ is 5.
		\item[(iv)] $X$ is 3-edge-colourable. Furthermore, in a proper 3-edge-colouring $f$ of $X$, $f((u_0)) = f((v_0))$ and $f((u_1)) = f((v_1)).$ 		 
	\end{enumerate}
\end{lemma}

\begin{remark} \label{remarkX} 
	{\rm Let $G$ be a graph which contains the semi-graph $X$ as in Figure \ref{x}. Let $C$ be a cycle in a 2-factor of $G$ which contains from $X$ the semi-edges $(v_0)$ and $(u_1)$. It is easy to establish that it is impossible for $C$, or any other cycle in the 2-factor, to contain both the other two semi-edges $(u_0)$ and $(v_1)$ from $X$. Given this, as well as Lemma \ref{lemmaX} (i) and \ref{lemmaX} (iii), the only possibility for $C$ is that it contains exactly one other vertex from $X$, so that the remaining vertices of $X$ may form a 5-cycle as part of the 2-factor of $G$. Thus, the remaining vertices of $X$ necessarily add one odd component to the 2-factor of $G$. The same applies for a cycle containing semi-edges $(v_0)$ and $(v_1)$, $(u_0)$ and $(v_1)$, or $(u_0)$ and $(u_1)$. Essentially, if a cycle from a 2-factor traverses through $X$, then the 2-factor must contain an odd component entirely contained in $X$ as well.}
\end{remark}

\begin{figure}[h!]
	\begin{center}
		\begin{tikzpicture}[every
			node/.style={draw,shape=circle,fill=black,text=white,scale=0.6},scale=0.4]
			\begin{scope} [xshift=0cm] 
				\foreach \i [count=\ii from 0] in {54,126,198,270,342}{
					\path (\i:15mm) node (p\ii) {};
				}
				\foreach \x in {0,...,4}{
					\tikzmath{
						integer \y;
						\y = mod(\x+2,5);
					}
					\draw (p\x) -- (p\y);
				}
				\node [fill=white,draw=none] (b0) at (-4,-4) {};
				\node [fill=white,draw=none] (b4) at (4,-4) {};
				\draw (b0) -- (b4);
				\node (b1) at ($(b0)!(p2)!(b4)$) {};
				\node (b2) at ($(b0)!(p3)!(b4)$) {};
				\node (b3) at ($(b0)!(p4)!(b4)$) {};
				\draw (p2) -- (b1);
				\draw (p3) -- (b2);
				\draw (p4) -- (b3);
				\path let \p1 = (p1) in node [fill=white,draw=none] (a0) at (-4,\y1) {};
				\path let \p1 = (p1) in node [fill=white,draw=none] (a1) at (4,\y1) {};
				\draw (p1) -- (a0);
				\draw (p0) -- (a1);
				
				\draw (b0) -- (b1);
				\draw (p1) -- (a0);
				
				\path [late options={name=p0,label={[color=black,scale=1.5]:$u_1$}}];
				\path [late options={name=b3,label={[color=black,scale=1.5]270:$v_1$}}];
				
			\end{scope}
			\begin{scope} [xshift=8cm] 
				\foreach \i [count=\ii from 0] in {54,126,198,270,342}{
					\path (\i:15mm) node (p\ii) {};
				}
				\foreach \x in {0,...,4}{
					\tikzmath{
						integer \y;
						\y = mod(\x+2,5);
					}
					\draw (p\x) -- (p\y);
				}
				\node [fill=white,draw=none] (b4) at (4,-4) {};
				\draw (b0) -- (b4);
				\node (b1) at ($(b0)!(p2)!(b4)$) {};
				\node (b2) at ($(b0)!(p3)!(b4)$) {};
				\node (b3) at ($(b0)!(p4)!(b4)$) {};
				\draw (p2) -- (b1);
				\draw (p3) -- (b2);
				\draw (p4) -- (b3);
				\path let \p1 = (p1) in node  (a0) at (-4,\y1) {};
				\path let \p1 = (p1) in node  [fill=white,draw=none] (a1) at (4,\y1) {};
				\draw (p1) -- (a0);
				\draw (p0) -- (a1);

				\path let \p1 = (p1) in node [fill=white,draw=none] (x0) at (-4,4) {};

				\draw (b3) -- (b4);
				\draw (a0) -- (x0);

				\path [late options={name=p1,label={[color=black,scale=1.5]90:$u_0$}}];
				\path [late options={name=b1,label={[color=black,scale=1.5]270:$v_0$}}];
				
				\path [late options={name=a0,label={[color=black,scale=1.5]270:$x$}}];
			\end{scope}
		\end{tikzpicture}
	\end{center}
	\caption{The semi-graph $Y$. }\label{y}
\end{figure}
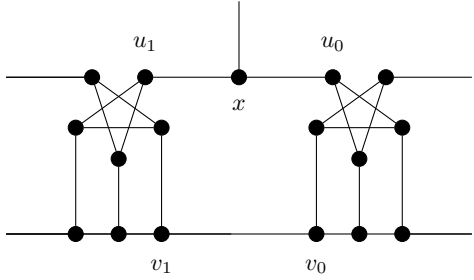

\begin{remark} \label{remarkY} 
	{\rm Consider the semi-graph $Y$ in Figure \ref{y}. Since any proper 3-edge-colouring would assign the same colours to $u_1x$ and $v_1v_0$ as per Lemma \ref{lemmaX}, and consequently to $u_0x$ as well, it is clear that $Y$ is not 3-edge-colourable. Indeed, $r(Y)=1$, where $x$ could be the only conflicting vertex in a 3-edge-colouring of $Y$. Also, $Y$ is Hamiltonian.}
\end{remark}

The following Figures \ref{fig_jerm2} and \ref{fig_jerm3} were presented as is in \cite{karabasetal}, although these constructions first appeared in \cite{kochol}. For our purposes, we utilise the properties of these graphs as they were proved in \cite{karabasetal}. We present brief explanations of the properties of these graphs in the following paragraph, taken from \cite{karabasetal}. For complete details and proofs, we advise the reader to consult the said paper.   


\begin{figure}[h!]
	\begin{center}
		\includegraphics[width=10cm]{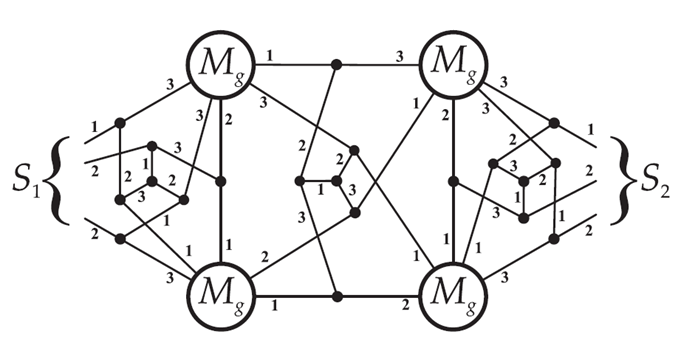}
	\end{center}
	\caption{A 3-edge-colourable cubic semi-graph $F_g$ with girth $g$ as presented in \cite{karabasetal}.} \label{fig_jerm2}
\end{figure}

\begin{figure} [h!]
	\begin{center}
		\includegraphics[width=10cm]{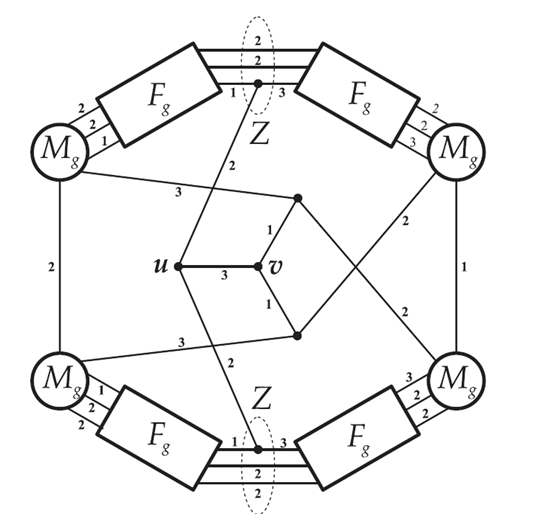}
	\end{center}
	\caption{A snark $G_g$ with resistance 2 and girth $g$ as presented in \cite{karabasetal}.} \label{fig_jerm3}
\end{figure}

The subgraph labelled $M_g$ present in the cubic semi-graph $F_g$ in Figure \ref{fig_jerm2} is itself a 3-edge-colourable cubic semi-graph with 5 semi-edges. The authors of \cite{karabasetal} created $M_g$ from a connected bipartite cubic graph $L_g$ of girth $g = 2r \geq 6$ by removing a path of length 2 (with semi-edges remaining). The semi-edges of each $M_g$ are joined arbitrarily in $F_g$, as mentioned in \cite{karabasetal} (see the fourth paragraph of the proof of Theorem 5.1 in \cite{karabasetal}). As the authors explained, the existence of such an $L_g$ is guaranteed as follows: Theorem 4.8 in \cite{nedela} guarantees the existence of an arc-transitive cubic graph $K$ of girth $g$. If $K$ is bipartite, let $L_g = K$. If $K$ is not bipartite, let $L_g$ be the bipartite double of $K$ (the direct product $K \times K_2$ with the complete graph $K_2$ on two vertices), which is connected, cubic, bipartite (and thus 3-edge-colourable), and has girth $g$. $M_g$ is obtained from $L_g$ by the removal of any path of length 2, with 5 semi-edges remaining. Given that the 5 edges joining to each $M_g$ are joining to 5 semi-edges remaining after the removal of a 2-path, and that there exists no cycles of $F_g$ which do not traverse through any of the $M_g$ subgraphs, we have that the girth, and also semi-girth, of $F_g$ is not smaller than $g$. $F_g$ was also constructed in such a way that if we joined the three semi-edges in $S_1$ to a single vertex, and the three semi-edges in $S_2$ to a single vertex, the resulting graph is a snark with resistance 2. Therefore, any proper 3-edge-colouring of $F_g$ must be such that two semi-edges in $S_1$ are coloured $a$, and one coloured $b$, for distinct $a,b \in \{1,2,3\}$ not necessarily distinct. Consider now the snark $G_g$ in Figure \ref{fig_jerm3}. Note that here as well the order in which the semi-edges of $M_g$ (those not contained in any $F_g$) are joined, is irrelevant with regards to whether $G_g$ is a snark (see the last paragraph of the proof of Theorem 5.1 in \cite{karabasetal} where the authors state the snark $G_g$ is not uniquely determined for this reason). Any cycle in $G$ must traverse through an instance of $M_g$, so we have that the girth of $G_g$ is also not smaller than $g$. The defect of $G_g$ is thus greater than $\lceil \frac{g}{2} \rceil$.

We now construct a cubic semi-graph $Z_g$ by removing any one the instances of $M_g$ from $G_g$ which are distance two from vertex $v$ in Figure \ref{fig_jerm3}, with 5 semi-edges remaining. 
If $Z_g$ is 3-edge-colourable, then we can properly colour $Z_g$ such that the 5 semi-edges are coloured in parity as per the Parity Lemma. That is, the colours of the 5 semi-edges would be $a,a,a,b,c$ where $a,b,c$ are distinct in $\{1,2,3\}$. In which case, given that the semi-edges of $M_g$ can be joined arbitrarily, we could add back a properly 3-edge-coloured $M_g$ and join semi-edges so that we have a proper 3-edge-colouring of $G_g$, a contradiction. Therefore, $Z_g$ is not 3-edge-colourable. It is easy to find a 3-edge-colouring of $Z_g$ with the only conflicting vertex being $u$ (and the 5 semi-edges could be coloured $a,a,b,b,c$ not necessarily distinct). Thus, $r(Z_g)=1$.  
Since the girth of $G_g$ is not smaller than $g$, it is also immediately clear that the semi-girth of $Z_g$ is not smaller than $g$. As such, the defect of $Z_g$ is thus greater than $\lceil \frac{g}{2} \rceil$.   

			
			
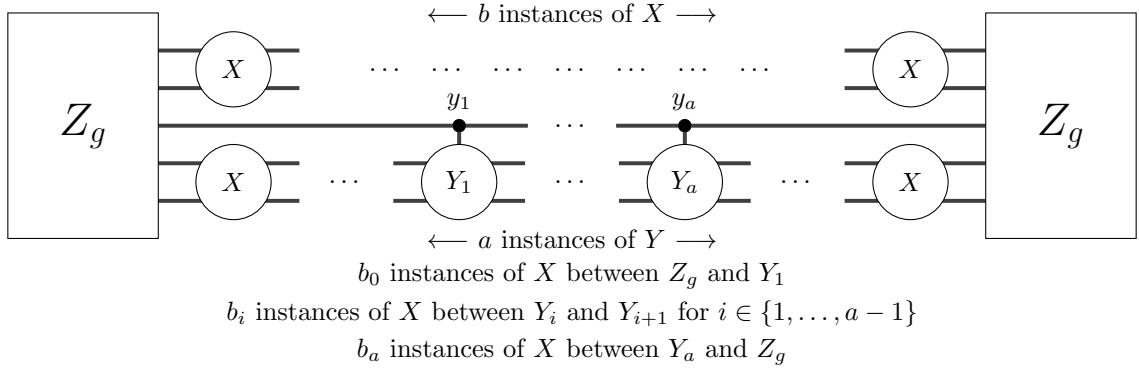
\begin{figure}[h!] 
	\begin{center}
		\begin{tikzpicture}[every node/.style={draw,shape=rectangle,fill=black,text=white},scale=1]
			\node[rectangle,draw,fill=white,text=black,minimum width = 2cm, minimum height = 3cm] (zg1) at (0,0) {\LARGE $Z_g$};	
			\node[rectangle,draw,fill=white,text=black,minimum width = 2cm, minimum height = 3cm] (zg2) at (13,0) {\LARGE $Z_g$};
			
			\node[fill=none,draw=none] (zg11) at (0,1) {};
			\node[fill=none,draw=none] (zg12) at (0,0.5) {};
			\node[fill=none,draw=none] (zg13) at (0,0) {};
			\node[fill=none,draw=none] (zg14) at (0,-0.5) {};
			\node[fill=none,draw=none] (zg15) at (0,-1) {};
			
			\node[fill=none,draw=none] (zg21) at (13,1) {};
			\node[fill=none,draw=none] (zg22) at (13,0.5) {};
			\node[fill=none,draw=none] (zg23) at (13,0) {};
			\node[fill=none,draw=none] (zg24) at (13,-0.5) {};
			\node[fill=none,draw=none] (zg25) at (13,-1) {};

			\node[circle,draw,fill=white,text=black,minimum width = 1cm, minimum height = 1cm] (x11) at (2,0.75) {$X$};
			\node[circle,draw,fill=white,text=black,minimum width = 1cm, minimum height = 1cm] (x12) at (2,-0.75) {$X$};
						
			\node[circle,draw,fill=white,text=black,minimum width = 1cm, minimum height = 1cm] (x21) at (11,0.75) {$X$};
			\node[circle,draw,fill=white,text=black,minimum width = 1cm, minimum height = 1cm] (x22) at (11,-0.75) {$X$};
			
			\node[circle,draw=none,fill=none,text=black] (text1) at (6.5,1.5) {$\longleftarrow$ $b$ instances of $X$ $\longrightarrow$};
			\node[circle,draw=none,fill=none,text=black] (text1) at (6.5,0.75) {$\dots$ \ \ $\dots$ \ \ $\dots$ \ \ $\dots$ \ \ $\dots$ \ \ $\dots$ \ \ $\dots$};

			\node[fill=black,draw,shape=circle,scale=0.5,label={[color=black]above:$y_1$}] (y1) at (5,0) {};
			\node[fill=black,draw,shape=circle,scale=0.5,label={[color=black]above:$y_a$}] (y2) at (8,0) {};
			
			\node[circle,draw,fill=white,text=black,minimum width = 1cm, minimum height = 1cm] (y11) at (5,-0.75) {$Y_1$};
			\node[circle,draw,fill=white,text=black,minimum width = 1cm, minimum height = 1cm] (y22) at (8,-0.75) {$Y_a$};
			
			\Edge[,bend=0](y1)(y11)
			\Edge[,bend=0](y2)(y22)	
		
			\Edge[,bend=0](zg13)(y1)
			\Edge[,bend=0](zg23)(y2)
		
			\node[fill=none,draw=none,shape=circle,scale=0.5] (w1) at (6,0) {};
			\node[fill=none,draw=none,shape=circle,scale=0.5] (w2) at (7,0) {};
			\node[fill=none,draw=none,text=black,shape=circle,scale=1] (w3) at (6.5,0) {$\dots$};
			\node[fill=none,draw=none,text=black,shape=circle,scale=1] (w4) at (6.5,-0.75) {$\dots$};
			
			\node[fill=none,draw=none,text=black,shape=circle,scale=1] (w5) at (3.5,-0.75) {$\dots$};
			\node[fill=none,draw=none,text=black,shape=circle,scale=1] (w6) at (9.5,-0.75) {$\dots$};
		
			\Edge[,bend=0](y1)(w1)
			\Edge[,bend=0](y2)(w2)
			
			\node[rectangle,draw=none,fill=none,text=black] (text2) at (6.5,-1.5) {$\longleftarrow$ $a$ instances of $Y$ $\longrightarrow$};
			\node[rectangle,draw=none,fill=none,text=black] (text2) at (6.5,-2) {$b_0$ instances of $X$ between $Z_g$ and $Y_1$}; 
			\node[rectangle,draw=none,fill=none,text=black] (text2) at (6.5,-2.5) {$b_i$ instances of $X$ between $Y_i$ and $Y_{i+1}$ for $i \in \{1,\dots,a-1\}$};
			\node[rectangle,draw=none,fill=none,text=black] (text2) at (6.5,-3) {$b_a$ instances of $X$ between $Y_a$ and $Z_g$};

			\node[fill=none,draw=none] (v1) at (3,1) {};
			\node[fill=none,draw=none] (v2) at (3,0.5) {};
			\node[fill=none,draw=none] (v4) at (3,-0.5) {};
			\node[fill=none,draw=none] (v5) at (3,-1) {};
			
			\node[fill=none,draw=none] (u1) at (10,1) {};
			\node[fill=none,draw=none] (u2) at (10,0.5) {};
			\node[fill=none,draw=none] (u4) at (10,-0.5) {};
			\node[fill=none,draw=none] (u5) at (10,-1) {};

			\node[fill=none,draw=none] (y1x) at (4,-1) {};
			\node[fill=none,draw=none] (y1y) at (4,-0.5) {};
			\node[fill=none,draw=none] (y2y) at (6,-0.5) {};
			\node[fill=none,draw=none] (y2x) at (6,-1) {};
			\node[fill=none,draw=none] (y3x) at (7,-1) {};
			\node[fill=none,draw=none] (y3y) at (7,-0.5) {};
			\node[fill=none,draw=none] (y4y) at (9,-0.5) {};
			\node[fill=none,draw=none] (y4x) at (9,-1) {};
		
			\Edge[,bend=0](y1x)(y2x)
			\Edge[,bend=0](y1y)(y2y)
			\Edge[,bend=0](y3x)(y4x)
			\Edge[,bend=0](y3y)(y4y)

			\Edge[,bend=0](zg11)(v1)
			\Edge[,bend=0](zg12)(v2)
			\Edge[,bend=0](zg14)(v4)
			\Edge[,bend=0](zg15)(v5)
		
			\Edge[,bend=0](zg21)(u1)
			\Edge[,bend=0](zg22)(u2)
			\Edge[,bend=0](zg24)(u4)
			\Edge[,bend=0](zg25)(u5)	
		
	\end{tikzpicture}
	\end{center}
	\caption{Cubic graph $G^a_g(B)$ where $B = \{b,b_0,\dots,b_a\}$.}\label{fig_Gabc}
\end{figure}

To prove the main result of this paper we utilise instances of semi-graph $X$, semi-graph $Y$ and semi-graph $Z_g$ to construct an infinite class of graphs with the desired properties. See Figure \ref{fig_Gabc}. 

\begin{theorem}
	For any integers $r \geq 3$, even $\omega \geq r$, and $d \geq \frac{3}{2}\omega$, there exists a non-trivial snark $G$ with $r(G)=r$, $\omega(G)=\omega$ and $\mu_3(G) \geq d$.
	\begin{proof}
		Let $r,\omega, d$ be positive integers such that $r \geq 3$, $\omega$ is even and $\omega \geq r$, and $d \geq \frac{3}{2}\omega$.  
		
		Consider the graph $G^a_g(B)$ as in Figure \ref{fig_Gabc}. Note that in the construction of this graph, the semi-edges of each instance of $X$ (including those in each instance of $Y_i$) are joined as per Figure \ref{x} and Figure \ref{y}. That is, semi-edge $(u_0)$ and $(v_0)$ on one side, and semi-edges $(u_1)$ and $(v_1)$ on the other side. Given that $r(Z_g)=1$ and $r(Y)=1$, it is clear that $r(G^a_g(B)) \geq a + 2$ whichever way the semi-edges of $Z_g$ are joined. Given: that any 3-edge-colouring of $Z_g$ with $r(Z_g)=1$ conflict will have semi-edges coloured $a,a,b,b,c$ not necessarily distinct; Lemma \ref{lemmaX} (iv); and Remark \ref{remarkY}, we have that it is possible to join $Z_g$ in $G^a_g(B)$ and colour edges with 3 colours in such a way that there is just one conflict in each $Z_g$ and the only other conflicting vertices are each $x$ in each $Y_i$ (as per Figure \ref{y}). Therefore, setting $a = r-2$, we can construct $G^a_g(B)$ such that $r(G^a_g(B))=r$.
		
		Consider some vertex $y_i$ incident to some $Y_i$ in $G^a_g(B)$. Let $F$ be a 2-factor of $G^a_g(B)$. The vertex $y_i$ is contained in some cycle $C$ in $F$. This cycle $C$ traverses to $y_{i-1}$ (or $Z_g$) on one side, and/or to $y_{i+1}$ (or $Z_g$) on the other side. It then necessarily has to traverse through: the $b$ instances of $X$ joining the two $Z_g's$; and/or the $b_{i-1}$ instances of $X$ joining $Y_{i-1}$ to $Y_i$; and/or the $b_{i}$ instances of $X$ joining $Y_{i}$ to $Y_{i+1}$. This implies that there must be at least one cycle in any 2-factor of $G^a_g(B)$ which has to traverse through at least $\min\{b,b_0,\dots,b_a\}$ instances of $X$. By Remark \ref{remarkX}, we then have that $\omega(G^a_g(B)) \geq b$. 
		
		Now, say $F$ is a 2-factor of $G_g$ with 2 odd components. $F$ must contain vertex $v$ and must therefore contain two or four joining edges to some instance of $M_g$. One of the two odd components in the 2-factor of $G_g$ must necessarily be disjoint from this instance of $M_g$ (otherwise, we could take $F$ as a 2-factor of $Z_g$, properly colour the edges in $F$ with colours 1 and 2 alternately, and properly colour the remaining edges with colour 3, which is a contradiction since $Z_g$ is not 3-edge-colourable). If we allow such an instance of $M_g$ to be the one removed to get semi-graph $Z_g$ from $G_g$, then we see that allowing the instance of $Z_g$ in $G^a_g(B)$ to inherit the edges of $F$, we may either close the components of $F$ in $G^a_g(B)$ by traversing through a Hamiltonian path of $X$, or traversing through vertex $y_1$ (or $y_a$) and returning to $Z_g$ by traversing through multiple instances of $X$. Either way, the instance of $Z_g$ adds only one odd component to a 2-factor of $G^a_g(B)$ with minimum odd components. Each instance of $Y$ also adds one odd component to a 2-factor of $G^a_g(B)$ with minimum odd components. To ensure that $\omega(G^a_g(B))=\omega$, we may simply adjust the values of $b_0,\dots,b_a$ and $b$ accordingly, so that the minimum number of odd components in any 2-factor equals exactly $\omega$. For example, we could set each of $b_0,b_1,\dots,b_a$ to be extremely large, and have $b$ be appropriately smaller than $\min\{b_0,\dots,b_a\}$ so that a 2-factor with minimum number of odd components necessarily has just one odd component in each $Z_g$, just one odd component in each the $b$ instances of $X$ joining the two $Z_g's$, and just one odd component in each $Y_i$.  
		
		Since $Z_g$ is a semi-subgraph of $G^a_g(B)$, and the defect of $Z_g$ is greater than $\lceil \frac{g}{2} \rceil$, we have by Lemma \ref{lem_semigirth} that $\mu_3(G^a_g(B)) \geq \lceil \frac{g}{2} \rceil$. We may simply choose $g$ to be large enough such that $\mu_3(G) \geq d$.
		
		$X$ has a girth 5, and it is clear that there are no cycles of size less than 5. It is also clear by inspection that there are no cyclic cuts of $G^a_g(B)$ of order less than 4. $G^a_g(B)$ is therefore a non-trivial snark.
	\end{proof}
\end{theorem}

\end{document}